\documentclass{amsart}
\usepackage{amssymb}
\usepackage{amsmath}
\usepackage{amsfonts}

\newtheorem{theorem}{Theorem}
\theoremstyle{plain}

\newtheorem{corollary}{Corollary}

\newtheorem{example}{Example}

\newtheorem{lemma}{Lemma}

\newtheorem{proposition}{Proposition}
\newtheorem{remark}{Remark}

\numberwithin{equation}{section}
\input{tcilatex}

\begin{document}
\title[Existence of Global Solutions]{Global existence of solutions for
Gierer-Meinhardt\ system with three equations}
\author{Abdelmalek Salem}
\curraddr{Department of Mathematics, University of Tebessa 12002 Algeria}
\email{a.salem@gawab.com}
\author{ Louafi Hichem}
\curraddr{Department of Mathematics , University of Annaba 23000 Algeria}
\email{hichemlouafi@gmail.com}
\author{Youkana Amar}
\curraddr{Department of Mathematics , University of Batna 05000 Algeria}
\email{youkana\_amar@yahoo.fr}
\thanks{}
\date{}
\dedicatory{}
\thanks{}

\begin{abstract}
This paper deals with an Gierer-Meinhardt model, with three substances,
formed Reaction-Diffusion system with fractional reaction. To prove global
existence for solutions of this system presents difficulties at the
boundednees of fractionar term. The purpose of this paper is to prove the
existence of a global solution using a boundary functionel. Our technique is
based on the construction of Lyapunov functionel.
\end{abstract}

\maketitle

\section{\textbf{Introduction}}

In recent years, systems of Reaction-Diffusion equations have received a
great deal of attention, motivated by their widespread occurrence in models
of chemical and biological phenomena. These systems are divided into
celebrated classes; systems\ with dissipation of mass and systems of
Gierer-Meinhardt. In this paper we deal with this last.

In the study of the various topics from plant developmental; Meinhardt, Koch
and Bernasconi\ \cite{Meinhardt3} proposed Activator-Inhibitor models (an
example is given in section 5) to describe a theory of biological pattern
formation in plants (\textit{Phyllotaxis}).

We assume a Reaction-diffusion system with three components:%
\begin{equation}
\left\{ 
\begin{array}{l}
\frac{\partial u}{\partial t}-a_{1}\Delta u=f(u,v,w)=\sigma -b_{1}u+\frac{%
u^{p_{1}}}{v^{q_{1}}(w^{r_{1}}+c)} \\ 
\frac{\partial v}{\partial t}-a_{2}\Delta v=g(u,v,w)=-b_{2}v+\frac{u^{p_{2}}%
}{v^{q_{2}}w^{r_{2}}} \\ 
\frac{\partial w}{\partial t}-a_{3}\Delta w=h(u,v,w)=-b_{3}w+\frac{u^{p_{3}}%
}{v^{q_{3}}w^{r_{3}}}%
\end{array}%
\right. \ x\in \Omega ,t>0\text{,}  \tag*{(1.1)}
\end{equation}%
with Neummann boundary conditions%
\begin{equation}
\frac{\partial u}{\partial \eta }=\frac{\partial v}{\partial \eta }=\frac{%
\partial w}{\partial \eta }=0\text{\ \ \ \ \ on }\partial \Omega \times
\left\{ t>0\right\} ,  \tag*{(1.2)}
\end{equation}%
and the initial data%
\begin{equation}
\left\{ 
\begin{array}{c}
u(0,x)=\varphi _{1}(x)>0 \\ 
v(0,x)=\varphi _{2}(x)>0 \\ 
w(0,x)=\varphi _{3}(x)>0%
\end{array}%
\right. \text{\ on}\;\Omega ,  \tag*{(1.3)}
\end{equation}%
and $\varphi _{i}\in C\left( \overline{\Omega }\right) $ for all $i=1,2,3.$

Here $\Omega $ is an open bounded domain of class $\mathbb{C}^{1}$ in $%
\mathbb{R}^{N}$, with boundary $\partial \Omega \;$and $\dfrac{\partial }{%
\partial \eta }$ denotes the outward normal derivative on $\partial \Omega
.\ $

$c,,p_{i},$ $q_{i}$, $r_{i}$: are non negative with $\sigma ,b_{i},$ $a_{i}$ 
$>0$, indexes for all $i=1,2,3.$%
\begin{equation}
0<p_{1}-1<\max \left\{ p_{2}\min \left( \frac{q_{1}}{q_{2}+1},\frac{r_{1}}{%
r_{2}},1\right) ,\text{ \ }p_{3}\min \left( \frac{r_{1}}{r_{3}+1},\frac{q_{1}%
}{q_{3}},1\right) \right\} .  \tag*{(1.4)}
\end{equation}

Put $A_{ij}=\frac{a_{i}+a_{j}}{2\sqrt{a_{i}a_{j}}}$ for all $i,j=1,2,3$. Let 
$\alpha ,\beta $ and $\gamma $ be positive constants such that~where

\begin{eqnarray}
\alpha &>&2\max \left\{ 1,\frac{b_{2}+b_{3}}{b_{1}}\right\} , 
\TCItag*{(1.5)} \\
\frac{1}{\beta } &>&2A_{12}^{2},  \TCItag*{(1.6)}
\end{eqnarray}

and%
\begin{equation}
\left( \frac{1}{2\beta }-A_{12}^{2}\right) \left( \frac{1}{2\gamma }%
-A_{13}^{2}\right) >\left( \frac{\alpha -1}{\alpha }A_{23}-A_{12}A_{13}%
\right) ^{2}.  \tag*{(1.7)}
\end{equation}

The main result of the paper reads as follows:

\begin{theorem}
Suppose that the functions $f,g$ and $h$ are satisfing condition (1.4). Let $%
\left( u\left( t,.\right) ,v\left( t,.\right) ,w\left( t,.\right) \right) $
be a solution of (1.1)-(1.3) and let:%
\begin{equation}
L(t)=\int_{\Omega }\frac{u^{\alpha }\left( t,x\right) }{v^{\beta }\left(
t,x\right) w^{\gamma }\left( t,x\right) }dx.  \tag*{(1.8)}
\end{equation}%
Then the functional $L$ is uniformly bounded on the interval $[0,T^{\ast
}],T^{\ast }<T_{\max }$.\newline
Where $T_{\max }$ $\left( \left\Vert u_{0}\right\Vert _{\infty },\left\Vert
v_{0}\right\Vert _{\infty },\left\Vert w_{0}\right\Vert _{\infty }\right) $
denotes the eventual blow-up time.
\end{theorem}

\begin{corollary}
Under the assumptions of theorem 1 all solutions of (1.1)-(1.3) with
positive initial data in $C\left( \overline{\Omega }\right) $ are global. If
in addition $b_{1}$, $b_{2}$, $b_{3}$, $\sigma >0,$ then $(u,v,w)$ are
uniformly bounded in $\overline{\Omega }\times \left[ 0,\infty \right) .$
\end{corollary}

\section{\textbf{Previous Results}}

The usual norms in spaces $L^{p}(\Omega )$, $L^{\infty }(\Omega )$ and $C(%
\overline{\Omega })$ are denoted respectively by:

\begin{eqnarray}
\left\Vert u\right\Vert _{p}^{p} &=&\frac{1}{\left\vert \Omega \right\vert }%
\int_{\Omega }\left\vert u(x)\right\vert ^{p}dx;  \TCItag*{(2.1a)} \\
\left\Vert u\right\Vert _{\infty } &=&\underset{x\in \Omega }{\max }%
\left\vert u(x)\right\vert ,  \TCItag*{(2.1b)} \\
\left\Vert u\right\Vert _{\mathbb{C}(\overline{\Omega })} &=&\underset{x\in 
\overline{\Omega }}{\max }\left\vert u(x)\right\vert .  \TCItag*{(2.1c)}
\end{eqnarray}

In 1972, following an ingenious idea of Turing.A \cite{Turing}, Gierer.A and
Meinhardt. H \cite{Meinhardt3} proposed a mathematical model for pattern
formations of spatial tissue structures of hydra in morphogenesis, a
biological phenomenon discovered by A.Trembley in 1744 \cite{Abraham}.It is
a system of reaction-diffusion equations of the form:%
\begin{equation}
\left\{ 
\begin{array}{l}
\frac{\partial u}{\partial t}-a_{1}\Delta u=\sigma -\mu u+\frac{u^{p}}{v^{q}}
\\ 
\frac{\partial v}{\partial t}-a_{2}\Delta v=-\nu v+\frac{u^{r}}{v^{s}}%
\end{array}%
\right. \text{ \ \ \ \ for all }x\in \Omega \text{\ },t>0  \tag*{(2.2)}
\end{equation}%
with Neummann boundary conditions%
\begin{equation}
\frac{\partial u}{\partial \eta }=0,\text{ and }\frac{\partial v}{\partial
\eta }=0,\text{ }\ \ \ \ x\in \partial \Omega ,t>0,  \tag*{(2.3)}
\end{equation}%
and initial conditions%
\begin{equation}
\left\{ 
\begin{array}{c}
u\left( x,0\right) =\varphi _{1}(x)>0 \\ 
v\left( x,0\right) =\varphi _{2}(x)>0%
\end{array}%
,\ \ \ \ \ \ \ x\in \Omega ,\right.  \tag*{(2.4)}
\end{equation}%
where $\Omega \subset 
\mathbb{R}
^{N}$ is a bounded domain with smooth boundary $\partial \Omega $, $%
a_{1},a_{2}>0,\mu ,\nu ,\sigma >0,p$, $q$, $r$ and $s$ are non negative
indexes with$\ p>1.$

Global existence of solutions in $\left( 0,\infty \right) $ is proved by
Rothe in 1984 \cite{Rothe} with special cases: $p=2,\ q=1,\ r=2,\ s=0$ and $%
N=3.$The Rothe's method cannot be applied (at least directly) to the general 
$p,q,r,s$. It is desirable to consider the $p,q,r,s$ originally proposed by
Gierer-Meinhardt. Wu and Li \cite{Wu} obtained the same results for the
problem (2.1)-(2.3) so long as $u,v^{-1}$ and $\sigma $\ are suitably small. 
\textbf{\ }Mingde, Shaohua and Yuchun \cite{Mingde}\ show that the solutions
of this problem are bounded all the time for each pair of initial values if 
\begin{eqnarray}
\frac{p-1}{r} &<&\frac{q}{s+1},  \TCItag*{(2.5a)} \\
\frac{p-1}{r} &<&1.  \TCItag*{(2.5b)}
\end{eqnarray}

Masuda. K and Takahashi. K \cite{Masuda} we consider a more general system
for $\left( u,v\right) :$ 
\begin{equation}
\left\{ 
\begin{array}{l}
\frac{\partial u}{\partial t}-a_{1}\Delta u=\sigma _{1}\left( x\right) -\mu
u+\rho _{1}\left( x,u\right) \frac{u^{p}}{v^{q}} \\ 
\frac{\partial v}{\partial t}-a_{2}\Delta v=\sigma _{2}\left( x\right) -\nu
v+\rho _{2}\left( x,u\right) \frac{u^{r}}{v^{s}}%
\end{array}%
\right.  \tag*{(2.6)}
\end{equation}%
with $\sigma _{1},\sigma _{2}\in C^{1}\left( \overline{\Omega }\right) ,$ $%
\sigma _{1}\geq 0,\sigma _{2}\geq 0,\rho _{1},\rho _{2}\in C^{1}\left( 
\overline{\Omega }\times \overline{%
\mathbb{R}
}_{+}^{2}\right) \cap L^{\infty }\left( \overline{\Omega }\times \overline{%
\mathbb{R}
}_{+}^{2}\right) $ satisfying $\rho _{1}\geq 0,\rho _{2}>0$ and $p$, $q$, $r$
,$s$ are nonnegative constants satisfying (2.5a). Obviously, (2.4) system is
a special case of (2.6) system. In 1987, Masuda. K and Takahashi. K \cite%
{Masuda} extended the result to $\frac{p-1}{r}<\frac{2}{N+2}$ under the
condition $\sigma _{1}>0.$In 2006 Jiang.H \cite{jiang} under the conditions
(2.5a) - (2.5b) , \ $\varphi _{1},\varphi _{2}\in W^{2,l}\left( \Omega
\right) ,l>\max \left\{ N,2\right\} ,$ $\frac{\partial \varphi _{1}}{%
\partial \eta }=\frac{\partial \varphi _{2}}{\partial \eta }=0$ on $\partial
\Omega $ and $\varphi _{1}\geq 0,\varphi _{2}>0$ in $\overline{\Omega }.$%
Then (2.6) system has a unique nonnegative global solution $(u,v)$
satisfying (2.3)-(2.4).

\section{\textbf{Preliminary Observations}}

It is well-known that to prove global existence of solutions to $(1.1)-(1.3)$
(see Henry \cite{Henry}), it suffices to derive a uniform estimate of $%
\left\Vert f\left( u,v,w\right) \right\Vert _{p},$ $\left\Vert g\left(
u,v,w\right) \right\Vert _{p}$ and $\left\Vert h\left( u,v,w\right)
\right\Vert _{p}$on $[0;T_{\max }[$ in the space $L^{p}(\Omega )$ for some $%
p>N/2.$Our aim is to construct polynomial Lyapunov functionals allowing us
to obtain $L^{p}-$ bounds on $u;v$ and $w$ that lead to global existence.
Since the functions $f,g$ and $h$ are continuously differentiable on $%
\mathbb{R}
_{+}^{3}$, then for any initial data in $C(\overline{\Omega })$, it is easy
to check directly their Lipschitz continuity on bounded subsets of the
domain of a fractional power of the operator 
\begin{equation}
O=-\left( 
\begin{array}{ccc}
a_{1}\Delta & 0 & 0 \\ 
0 & a_{2}\Delta & 0 \\ 
0 & 0 & a_{3}\Delta%
\end{array}%
\right) .  \tag*{(3.1)}
\end{equation}%
Under these assumptions, the following local existence result is well known
(see Friedman \cite{Friedman} and Pazy \cite{Pazy}).

\begin{proposition}
The system $(1.1)-(1.3)$ admits a unique, classical solution $(u;v;w)$ on$\
(0,T_{\max }[\times \Omega $.%
\begin{equation}
\text{If }T_{\max }<\infty \text{ then }\underset{t\nearrow T_{\max }}{\lim }%
\left( \left\Vert u\left( t,.\right) \right\Vert _{\infty }+\left\Vert
v\left( t,.\right) \right\Vert _{\infty }+\left\Vert w\left( t,.\right)
\right\Vert _{\infty }\right) =\infty \text{.}  \tag*{(3.2)}
\end{equation}
\end{proposition}

\begin{remark}
This proposition seems to be well-known (Dan Henry \cite{Henry}).
Neverthless we could not find it in the literature in the form stated here
and stated in the book of Franz Rothe (\cite{Rothe} pp: 111-118 with proof).
Usually the explosion property $(3.2)$ is only stated for some norm
involving smoothness, but not the $L_{\infty }-$norm.
\end{remark}

\section{\textbf{Proofs}}

For the proof of theorem 1, we need a preparatory Lemmas, which are proved
in the appendix.

\begin{lemma}
Assume that $p,$ $q,$ $r,$ $s,$ $m,$ and $n$ satisfying%
\begin{equation*}
\frac{p-1}{r}<\min \left( \frac{q}{s+1},\frac{m}{n},1\right) .
\end{equation*}%
For all $h$, $l$, $\alpha $, $\beta $, $\gamma >0$, there exist $C=C\left(
h,l,\alpha ,\beta ,\gamma \right) >0$ and $\theta =\theta \left( \alpha
\right) \in \left( 0,1\right) $, such that%
\begin{equation}
\alpha \frac{x^{p-1+\alpha }}{y^{q+\beta }z^{m+\gamma }}\leq \beta \frac{%
x^{r+\alpha }}{y^{s+1+\beta }z^{n+\gamma }}+C\left( \frac{x^{\alpha }}{%
y^{\beta }z^{\gamma }}\right) ^{\theta },\ \ \ \ x\geq 0,y\geq h,z\geq l~ 
\tag*{(4.1)}
\end{equation}
\end{lemma}

\begin{lemma}
Let $\mu ,T>0$ and $f_{j}=f_{j}(t)$ be a non-negative integrable function on 
$[0,T)$and $0<\theta _{j}<1$ \ ($j=1,...,J$). Let $W=W(t)$be a positive
function on $[0,T)$ satisfying the differential inequality%
\begin{equation}
\frac{dW\left( t\right) }{dt}\leq -\mu W\left( t\right) +\underset{j=1}{%
\overset{J}{\sum }}f_{j}(t)W^{\theta _{j}}\left( t\right) ,\text{ \ \ \ }%
0\leq t<T.  \tag*{(4.2)}
\end{equation}%
Then, we obtain that%
\begin{equation}
W\left( t\right) \leq \kappa ,\text{ \ \ \ }0\leq t<T,  \tag*{(4.3)}
\end{equation}%
where $\kappa $ is the maximal root of the following algebraic equation:%
\begin{equation}
x-\underset{j=1}{\overset{J}{\sum }}\left( \underset{0<t<T}{\sup }%
\int_{0}^{t}e^{-\mu \left( t-\xi \right) }f_{j}(\xi )d\xi \right) x^{\theta
_{j}}=W\left( 0\right) .  \tag*{(4.4)}
\end{equation}%
Moreover, if $T=+\infty ,$ then 
\begin{equation*}
\underset{t\nearrow \infty }{\lim \sup }W\left( t\right) \leq \kappa
_{\infty },
\end{equation*}%
where $\kappa _{\infty }$ is the maximal root of the following algebraic
equation: 
\begin{equation*}
x-\underset{j=1}{\overset{J}{\sum }}\left( \underset{t\nearrow \infty }{\lim
\sup }\int_{0}^{t}e^{-\mu \left( t-\xi \right) }f_{j}(\xi )d\xi \right)
x^{\theta _{j}}=0.
\end{equation*}
\end{lemma}

\begin{lemma}
Let $\left( u\left( t,.\right) ,v\left( t,.\right) ,w\left( t,.\right)
\right) $ be a solution of (1.1)-(1.3), then\ for any $(t,x)$\ in\ $%
(0,T_{\max }[\times \Omega $ we get 
\begin{equation}
\left\{ 
\begin{array}{c}
u(t,x)\geq e^{-b_{1}t}\min (\varphi _{1}(x))>0, \\ 
v(t,x)\geq e^{-b_{2}t}\min (\varphi _{2}(x))>0, \\ 
w(t,x)\geq e^{-b_{3}t}\min (\varphi _{3}(x))>0.%
\end{array}%
\right. \text{ \ }  \tag*{(4.5)}
\end{equation}
\end{lemma}

\begin{proof}[proof of theorem 1]
Differentiating $L\left( t\right) $ with respect to $t$ yields%
\begin{eqnarray*}
L^{\prime }\left( t\right) &=&\int_{\Omega }\frac{d}{dt}\left( \frac{%
u^{\alpha }}{v^{\beta }w^{\gamma }}\right) dx, \\
&=&\int_{\Omega }\left( \alpha \frac{u^{\alpha -1}}{v^{\beta }w^{\gamma }}%
\partial _{t}u-\beta \frac{u^{\alpha }}{v^{\beta +1}w^{\gamma }}\partial
_{t}v-\gamma \frac{u^{\alpha }}{v^{\beta }w^{\gamma +1}}\partial
_{t}w\right) dx,
\end{eqnarray*}%
replacing $\partial _{t}u,$ $\partial _{t}v$\ and $\partial _{t}w$ with its
values in (1.1), we get%
\begin{eqnarray*}
L^{\prime }\left( t\right) &=&\int_{\Omega }\left( 
\begin{array}{c}
a_{1}\alpha \frac{u^{\alpha -1}}{v^{\beta }w^{\gamma }}\Delta u-a_{2}\beta 
\frac{u^{\alpha }}{v^{\beta +1}w^{\gamma }}\Delta v-a_{3}\gamma \frac{%
u^{\alpha }}{v^{\beta }w^{\gamma +1}}\Delta w \\ 
-b_{1}\alpha \frac{u^{\alpha }}{v^{\beta }w^{\gamma }}+b_{2}\beta \frac{%
u^{\alpha }}{v^{\beta }w^{\gamma }}+b_{3}\gamma \frac{u^{\alpha }}{v^{\beta
}w^{\gamma }} \\ 
+\alpha \frac{u^{p_{1}+\alpha -1}}{v^{q_{1}+\beta }w^{\gamma }(w^{r_{1}}+c)}%
-\beta \frac{u^{p_{2}+\alpha }}{v^{q_{2}+\beta +1}w^{r_{2}+\gamma }}-\gamma 
\frac{u^{p_{3}+\alpha }}{v^{q_{3}+\beta }w^{r_{3}+\gamma +1}}+\sigma \alpha 
\frac{u^{\alpha -1}}{v^{\beta }w^{\gamma }}%
\end{array}%
\right) dx, \\
&=&I+J,
\end{eqnarray*}%
where $I$ contains laplacian terms and $J$ contains the other terms 
\begin{eqnarray*}
I &=&a_{1}\alpha \int_{\Omega }\frac{u^{\alpha -1}}{v^{\beta }w^{\gamma }}%
\Delta udx-a_{2}\beta \int_{\Omega }\frac{u^{\alpha }}{v^{\beta +1}w^{\gamma
}}\Delta vdx-a_{3}\gamma \int_{\Omega }\frac{u^{\alpha }}{v^{\beta
}w^{\gamma +1}}\Delta wdx, \\
J &=&\left( -b_{1}\alpha +b_{2}\beta +b_{3}\gamma \right) L\left( t\right) \\
&&+\alpha \int_{\Omega }\frac{u^{p_{1}+\alpha -1}}{v^{q_{1}+\beta
}w_{3}^{\gamma }(w^{r_{1}}+c)}dx-\beta \int_{\Omega }\frac{u^{p_{2}+\alpha }%
}{v^{q_{2}+\beta +1}w^{r_{2}+\gamma }}dx-\gamma \int_{\Omega }\frac{%
u^{p_{3}+\alpha }}{v^{q_{3}+\beta }w^{r_{3}+\gamma +1}}dx \\
&&+\sigma \alpha \int_{\Omega }\frac{u^{\alpha -1}}{v^{\beta }w^{\gamma }}dx.
\end{eqnarray*}%
\newline
Starting with estimation of\textbf{\ }$I$:

Using Green's formula for terms $\int_{\Omega }\frac{u^{\alpha -1}}{v^{\beta
}w^{\gamma }}\Delta udx$ , $\int_{\Omega }\frac{u^{\alpha }}{v^{\beta
+1}w^{\gamma }}\Delta vdx$ and $\int_{\Omega }\frac{u^{\alpha }}{v^{\beta
}w^{\gamma +1}}\Delta wdx$ we get 
\begin{eqnarray*}
I &=&\int_{\Omega }\left( 
\begin{array}{c}
-a_{1}\alpha \left( \alpha -1\right) \frac{u^{\alpha -2}}{v^{\beta
}w^{\gamma }}\left\vert \nabla u\right\vert ^{2}+a_{1}\alpha \beta \frac{%
u^{\alpha -1}}{v^{\beta +1}w^{\gamma }}\nabla u\nabla v+a_{1}\alpha \gamma 
\frac{u^{\alpha -1}}{v^{\beta }w^{\gamma +1}}\nabla u\nabla w \\ 
+a_{2}\beta \alpha \frac{u^{\alpha -1}}{v^{\beta +1}w^{\gamma }}\nabla
u\nabla v-a_{2}\beta \left( \beta +1\right) \frac{u^{\alpha }}{v^{\beta
+2}w^{\gamma }}\left\vert \nabla v\right\vert ^{2}-a_{2}\beta \gamma \frac{%
u^{\alpha }}{v^{\beta +1}w^{\gamma +1}}\nabla v\nabla w \\ 
+a_{3}\gamma \alpha \frac{u^{\alpha -1}}{v^{\beta }w^{\gamma +1}}\nabla
u\nabla w-a_{3}\gamma \beta \frac{u^{\alpha }}{v^{\beta +1}w^{\gamma +1}}%
\nabla v\nabla w-a_{3}\gamma \left( \gamma +1\right) \frac{u^{\alpha }}{%
v^{\beta }w^{\gamma +2}}\left\vert \nabla w\right\vert ^{2}%
\end{array}%
\right) dx, \\
&=&-\int_{\Omega }\left[ \frac{u^{\alpha -2}}{v^{\beta +2}w^{\gamma +2}}%
\left( QT\right) \cdot T\right] dx,
\end{eqnarray*}%
where%
\begin{equation*}
Q=\left( 
\begin{array}{ccc}
a_{1}\alpha \left( \alpha -1\right) & -\alpha \beta \frac{a_{1}+a_{2}}{2} & 
-\alpha \gamma \frac{a_{1}+a_{3}}{2} \\ 
-\alpha \beta \frac{a_{1}+a_{2}}{2} & a_{2}\beta \left( \beta +1\right) & 
\beta \gamma \frac{a_{2}+a_{3}}{2} \\ 
-\alpha \gamma \frac{a_{1}+a_{3}}{2} & \beta \gamma \frac{a_{2}+a_{3}}{2} & 
a_{3}\gamma \left( \gamma +1\right)%
\end{array}%
\right) .
\end{equation*}%
$Q$ is matrix of quadratic form compared to: $vw\nabla u,$ $uw\nabla v$ and $%
uv\nabla w$, which is written in the matric form:$T=\left( vw\nabla
u,uw\nabla v,uv\nabla w\right) ^{t}$.\newline
$Q$ is definite positive if, and only if all its principal successive
determinants are positive. To see this, we have:\newline
1. $\Delta _{1}=a_{1}\alpha \left( \alpha -1\right) >0$. Using (1.5), we get 
$\Delta _{1}>0$.\newline
2. $\Delta _{2}=\left\vert 
\begin{array}{cc}
a_{1}\alpha \left( \alpha -1\right) & -\alpha \beta \frac{a_{1}+a_{2}}{2} \\ 
-\alpha \beta \frac{a_{1}+a_{2}}{2} & a_{2}\beta \left( \beta +1\right)%
\end{array}%
\right\vert =\alpha ^{2}\beta ^{2}a_{1}a_{2}\left( \frac{\alpha -1}{\alpha }%
\frac{\beta +1}{\beta }-A_{12}^{2}\right) $. Using (1.5) and (1.6), we get $%
\Delta _{2}>0$.\newline
3. Using theorem 1 in S.abdelmalek and S. Kouachi \cite{abdelmalek1} we get\ 
$\left( \alpha -1\right) \Delta _{3}=\left( \alpha -1\right) \left\vert
Q\right\vert =\alpha \left( \alpha \gamma \beta \right)
^{2}a_{1}a_{2}a_{3}(\left( \frac{\alpha -1}{\alpha }\frac{\beta +1}{\beta }%
-A_{12}^{2}\right) \left( \frac{\alpha -1}{\alpha }\frac{\gamma +1}{\gamma }%
-A_{13}^{2}\right) -\left( \frac{\alpha -1}{\alpha }A_{23}-A_{12}A_{13}%
\right) ^{2})$.\ Using (1.5)-(1.7), we get $\Delta _{3}>0.$\newline
Consequently we have $I\leq 0,$ \ \ \ \ \ $~\forall ~(t,x)~\in ~\left[
0,T^{\ast }\right] \times \Omega $.\newline
Now we estimate\textbf{\ }$J:$%
\begin{eqnarray*}
J &=&\left( -b_{1}\alpha +b_{2}\beta +b_{3}\gamma \right) L\left( t\right) \\
&&+\alpha \int_{\Omega }\frac{u^{p_{1}+\alpha -1}}{v^{q_{1}+\beta }w^{\gamma
}(w^{r_{1}}+c)}dx-\beta \int_{\Omega }\frac{u^{p_{2}+\alpha }}{%
v^{q_{2}+\beta +1}w^{r_{2}+\gamma }}dx-\gamma \int_{\Omega }\frac{%
u^{p_{3}+\alpha }}{v^{q_{3}+\beta }w^{r_{3}+\gamma +1}}dx \\
&&+\sigma \alpha \int_{\Omega }\frac{u^{\alpha -1}}{v^{\beta }w^{\gamma }}dx.
\end{eqnarray*}%
According to the maximum principle, there exists $C_{0}$ dependant on $%
\varphi _{1}$, $\varphi _{2}$ and $\varphi _{3}$ such that $v,w\geq C_{0}>0,$
then we have%
\begin{equation*}
\frac{u^{\alpha -1}}{v^{\beta }w^{\gamma }}=\left( \frac{u^{\alpha }}{%
v^{\beta }w^{\gamma }}\right) ^{\frac{\alpha -1}{\alpha }}\left( \frac{1}{v}%
\right) ^{\frac{\beta }{\alpha }}\left( \frac{1}{w}\right) ^{\frac{\gamma }{%
\alpha }}\leq \left( \frac{u^{\alpha }}{v^{\beta }w^{\gamma }}\right) ^{%
\frac{\alpha -1}{\alpha }}\left( \frac{1}{C_{0}}\right) ^{\frac{\beta
+\gamma }{\alpha }},
\end{equation*}%
then%
\begin{equation*}
\frac{u^{\alpha -1}}{v^{\beta }w^{\gamma }}\leq C_{2}\left( \frac{u^{\alpha }%
}{v^{\beta }w^{\gamma }}\right) ^{\frac{\alpha -1}{\alpha }}\ \ \ \ \text{%
where }C_{2}=\left( \frac{1}{C_{0}}\right) ^{\frac{\beta +\gamma }{\alpha }},
\end{equation*}%
we have%
\begin{eqnarray*}
J &=&\left( -b_{1}\alpha +b_{2}\beta +b_{3}\gamma \right) L\left( t\right) \\
&&+\alpha \int_{\Omega }\frac{u^{p_{1}+\alpha -1}}{v^{q_{1}+\beta }w^{\gamma
}\left( w^{r_{1}}+c\right) }dx-\beta \int_{\Omega }\frac{u^{p_{2}+\alpha }}{%
v^{q_{2}+\beta +1}w^{r_{2}+\gamma }}dx-\gamma \int_{\Omega }\frac{%
u^{p_{3}+\alpha }}{v^{q_{3}+\beta }w^{r_{3}+\gamma +1}}dx \\
&&+\sigma \alpha \int_{\Omega }\frac{u^{\alpha -1}}{v^{\beta }w^{\gamma }}dx.
\end{eqnarray*}%
\newline
Using lemma\ 1, $\forall (t,x)\in \left[ 0,T^{\ast }\right] [\times \Omega $
we get%
\begin{equation}
\alpha \frac{u^{p_{1}+\alpha -1}}{v^{q_{1}+\beta }w^{\gamma }(w^{r_{1}}+c)}%
\leq \alpha \frac{u^{p_{1}+\alpha -1}}{v^{q_{1}+\beta }w^{\gamma +r_{1}}}%
\leq \beta \frac{u^{p_{2}+\alpha }}{v^{q_{2}+\beta +1}w^{r_{2}+\gamma }}%
+C\left( \frac{u^{\alpha }}{v^{\beta }w^{\gamma }}\right) ^{\theta } 
\tag*{(4.6)}
\end{equation}%
or%
\begin{equation}
\alpha \frac{u^{p_{1}+\alpha -1}}{v^{q_{1}+\beta }w^{\gamma +r_{1}}}\leq
\gamma \frac{u^{p_{3}+\alpha }}{w^{r_{3}+1+\gamma }v^{q_{3}+\beta }}+C\left( 
\frac{u^{\alpha }}{v^{\beta }w^{\gamma }}\right) ^{\theta }  \tag*{(4.7)}
\end{equation}

Using (4.6) or (4.7) then%
\begin{equation*}
J\leq \left( -b_{1}\alpha +b_{2}\beta +b_{3}\gamma \right) L\left( t\right)
+\int_{\Omega }C\left( \frac{u^{\alpha }}{v^{\beta }w^{\gamma }}\right)
^{\theta }dx+\alpha \sigma \int_{\Omega }C_{2}\left( \frac{u^{\alpha }}{%
v^{\beta }w^{\gamma }}\right) ^{\frac{\alpha -1}{\alpha }}dx.
\end{equation*}%
Applying\ H\"{o}lder's inequality, for all $t$\ \ in$\ \left[ 0,T^{\ast }%
\right] \ $we\ obtain 
\begin{equation*}
\int_{\Omega }C\left( \frac{u^{\alpha }}{v^{\beta }w^{\gamma }}\right)
^{\theta }dx\leq \left( \int_{\Omega }\left( \frac{u^{\alpha }}{v^{\beta
}w^{\gamma }}\right) dx\right) ^{\theta }\left( \int_{\Omega }C^{\frac{1}{%
1-\theta }}dx\right) ^{1-\theta },
\end{equation*}%
then%
\begin{equation*}
\int_{\Omega }C\left( \frac{u^{\alpha }}{v^{\beta }w^{\gamma }}\right)
^{\theta }dx\leq C_{3}L^{\theta }(t),\ \ \ \ \ \ \text{where }%
C_{3}=C\left\vert \Omega \right\vert ^{1-\theta }.
\end{equation*}%
We have%
\begin{equation*}
\int_{\Omega }C_{2}\left( \frac{u^{\alpha }}{v^{\beta }w^{\gamma }}\right) ^{%
\frac{\alpha -1}{\alpha }}dx\leq \left( \int_{\Omega }\left( \frac{u^{\alpha
}}{v^{\beta }w^{\gamma }}\right) dx\right) ^{\frac{\alpha -1}{\alpha }%
}\left( \int_{\Omega }\left( C_{2}\right) ^{\alpha }dx\right) ^{\frac{1}{%
\alpha }},
\end{equation*}%
then%
\begin{equation*}
\int_{\Omega }C_{2}\left( \frac{u^{\alpha }}{v^{\beta }w^{\gamma }}\right) ^{%
\frac{\alpha -1}{\alpha }}dx\leq C_{4}L^{\frac{\alpha -1}{\alpha }}\left(
t\right) \ \ \text{where}\ C_{4}=C_{2}\left\vert \Omega \right\vert ^{\frac{1%
}{\alpha }},
\end{equation*}%
we get%
\begin{equation*}
J\leq \left( -b_{1}\alpha +b_{2}\beta +b_{3}\gamma \right) L\left( t\right)
+C_{3}L^{\theta }\left( t\right) +\alpha \sigma C_{4}L^{\frac{\alpha -1}{%
\alpha }}\left( t\right) ,
\end{equation*}%
which implies%
\begin{equation*}
J\leq \left( -b_{1}\alpha +b_{2}\beta +b_{3}\gamma \right) L\left( t\right)
+C_{5}\left( L^{\theta }\left( t\right) +\alpha \sigma L^{\frac{\alpha -1}{%
\alpha }}\left( t\right) \right) .
\end{equation*}%
\newline
Thus under conditions (1.5),\ (1.6) and (1.7), we obtain%
\begin{equation*}
L^{\prime }(t)\leq \left( -b_{1}\alpha +b_{2}\beta +b_{3}\gamma \right)
L\left( t\right) +C_{5}\left( L^{\theta }\left( t\right) +\alpha \sigma L^{%
\frac{\alpha -1}{\alpha }}\left( t\right) \right) ,
\end{equation*}%
\newline
since $-b_{1}\alpha +b_{2}\beta +b_{3}\gamma <0$ and Using lemma 2 we deduce
that $L(t)$ is bounded\ on\ $(0,T_{\max }[$\ ie $L(t)\leq \kappa $, where $%
\kappa $ dependent on $\varphi _{1},$ $\varphi _{2}$ and $\varphi _{3}$.
\end{proof}

\begin{proof}[Proof of Corollary 1]
Since $L(t)$ is bounded~on $(0,T_{\max }[\ $and\ the functions $\frac{%
u^{p_{1}}}{v^{q_{1}}(w^{r_{1}}+c)}$ $\frac{u^{p_{2}}}{v^{q_{2}}w^{r_{2}}}$
and $\frac{u^{p_{3}}}{v^{q_{3}}w^{r_{3}}}$ are$\ $in$\ $ $L^{\infty }((0$, $%
T_{\max })$, $L^{m}(%
\Omega
))$ for all $m>\frac{N}{2},$ then as a consequence of the arguments in
Henry. D \cite{Henry} or Haraux. A and Kirane. M \cite{Haraux}~we conclude
the solution of the system (1-1)-(1-7) is global and uniformly bounded on $%
\Omega \times (0,+\infty ).$
\end{proof}

\section{Example}

In this section we\ will examine a particular activator-inhibitor model in
order to illusrate the applicability of corollary 1 and proposition 1.We
assume that all reactions take place in a bounded domain $\Omega $ with a
smooth boundary $\partial \Omega $.

\begin{example}
The model proposed by Meinhardt, Koch and Bernasconi\ \cite{Meinhardt3}\ to
describe a theory of biological pattern formation in plants (\textit{%
Phyllotaxis}), where $u$, $v$ and $w$ are the concentrations of three
substances; called activator ($u$) and inhibitors ($v$ and $w$) is:%
\begin{equation}
\left\{ 
\begin{array}{l}
\frac{\partial u}{\partial t}-a_{1}\frac{\partial ^{2}u}{\partial x^{2}}%
=-b_{1}u+\frac{a^{2}}{v(w+\kappa _{u})}+\sigma , \\ 
\frac{\partial v}{\partial t}-a_{2}\frac{\partial ^{2}v}{\partial x^{2}}%
=-b_{2}v+u^{2}, \\ 
\frac{\partial w}{\partial t}-a_{3}\frac{\partial ^{2}w}{\partial x^{2}}%
=-b_{3}w+u,%
\end{array}%
\right. \text{for all }x\in \Omega ,\ t>0.  \tag*{(5.1)}
\end{equation}
\end{example}

\begin{proposition}
Solutions of $(5.1)$ with boundary conditions $(1.2)$ and nonnegative
uniformly bounded initial data $(1.3)$ exist globally.
\end{proposition}

\proof
This model is a special case of our general model (1.1), where $%
p_{1}=2,q_{1}=1,r_{1}=1,p_{2}=2,q_{2}=0,r_{2}=0,p_{3}=1,q_{3}=0,r_{3}=0.$
These indexes realize the conditions of global existence: $\frac{p_{1}-1}{%
p_{2}}<\min \left( \frac{q_{1}}{q_{2}+1},\frac{r_{1}}{r_{2}},1\right) .$%
\endproof%

\begin{remark}
The system described by equations $(5.1)$ exhibits all the essential
features of phyllotaxis.
\end{remark}

\section{Appendix}

The purpose of this appendix is to prove lemma 1, lemma 2 and lemma 3 in
section 4 which we have used in the proof of theorem 1.

\begin{proof}[Proof of Lemma 1]
For all $~x\geq 0,y\geq h,$ $z\geq l~~\ $we have from the inequality (4.1)

\begin{equation}
\alpha \frac{x^{p-1}}{y^{q}z^{m}}\leq \beta \frac{x^{r}}{y^{s+1}z^{n}}%
+C\left( \frac{x^{\alpha }}{y^{\beta }z^{\gamma }}\right) ^{\theta -1}~~~ 
\tag*{(6.1)}
\end{equation}%
and we can write%
\begin{equation*}
\alpha \frac{x^{p-1}}{y^{q}z^{m}}=\alpha \beta ^{-\frac{p-1}{r}}\left( \beta 
\frac{x^{r}}{y^{s+1}z^{n}}\right) ^{\frac{p-1}{r}}y^{\frac{\left( s+1\right)
\left( p-1\right) }{r}-q}z^{\frac{n\left( p-1\right) }{r}-m}~~.
\end{equation*}%
\newline
For each $\epsilon $ realize: $0<\epsilon <\min \left( \frac{q}{s+1},\frac{m%
}{n},1\right) -\frac{p-1}{r}$%
\begin{equation*}
\alpha \frac{x^{p-1}}{y^{q}z^{m}}=\alpha \beta ^{-\frac{p-1}{r}}\left( \beta 
\frac{x^{r}}{y^{s+1}z^{n}}\right) ^{\frac{p-1}{r}+\epsilon }\left( \beta 
\frac{x^{r}}{y^{s+1}z^{n}}\right) ^{-\epsilon }v^{\frac{\left( s+1\right)
\left( p-1\right) }{r}-q}z^{\frac{n\left( p-1\right) }{r}-m}.
\end{equation*}%
\newline
Then also 
\begin{eqnarray}
\alpha \frac{x^{p-1}}{y^{q}z^{m}} &=&\alpha \left( \beta \right) ^{-\frac{p-1%
}{r}-\epsilon }\left( \beta \frac{x^{r}}{y^{s+1}z^{n}}\right) ^{\frac{p-1}{r}%
+\epsilon }\left( \frac{1}{x^{\alpha }}\right) ^{\frac{r\epsilon }{\alpha }%
}\left( y\right) ^{\frac{\left( s+1\right) \left( p-1\right) }{r}-q+\epsilon
\left( s+1\right) }z^{\frac{n\left( p-1\right) }{r}-m+\epsilon n},  \notag \\
&\leq &\alpha \left( \beta \right) ^{-\frac{p-1}{r}-\epsilon }\left( \beta 
\frac{x^{r}}{y^{s+1}z^{n}}\right) ^{\frac{p-1}{r}+\epsilon }\left( \frac{1}{%
x^{\alpha }}\right) ^{\frac{r\epsilon }{\alpha }}\left( h\right) ^{\frac{%
\left( s+1\right) \left( p-1\right) }{r}-q+\epsilon \left( s+1\right) }l^{%
\frac{n\left( p-1\right) }{r}-m+\epsilon n},  \notag \\
&\leq &\alpha \left( \beta \right) ^{-\frac{p-1}{r}-\epsilon }\left( \beta 
\frac{x^{r}}{y^{s+1}z^{n}}\right) ^{\frac{p-1}{r}+\epsilon }\left( \frac{1}{%
x^{\alpha }}\right) ^{\frac{r\epsilon }{\alpha }}\left( h\right) ^{\frac{%
\left( s+1\right) \left( p-1\right) }{r}-q+\epsilon \left( s+1\right) }\times
\notag \\
&&l^{\frac{n\left( p-1\right) }{r}-m+\epsilon n}\left( \frac{y}{h}\right) ^{%
\frac{\beta r\epsilon }{\alpha }}\left( \frac{z}{l}\right) ^{\frac{\gamma
r\epsilon }{\alpha }},  \notag \\
&\leq &C_{1}\left( \beta \frac{x^{r}}{y^{s+1}z^{n}}\right) ^{\frac{p-1}{r}%
+\epsilon }\left( \frac{y^{\beta }z^{\gamma }}{x^{\alpha }}\right) ^{\frac{%
r\epsilon }{\alpha }},  \TCItag*{(6.2)}
\end{eqnarray}%
\newline
where%
\begin{equation*}
C_{1}=\alpha \left( \beta \right) ^{-\frac{p-1}{r}-\epsilon }h^{^{\frac{%
(s+1)(p-1)}{r}-q+\epsilon (s+1)-\frac{\beta r\epsilon }{\alpha }}}l^{\frac{%
(n)(p-1)}{r}-m+\epsilon n-\frac{\gamma r\epsilon }{\alpha }}.
\end{equation*}%
Using Young's inequality for (6.2) with taking $C=C_{1}^{1+\frac{%
p-1+r\epsilon }{r-\left( p-1\right) -r\epsilon }}$ and $\theta =1-\frac{%
r\epsilon }{\alpha \left( 1-\frac{p-1}{r}-\epsilon \right) }$ where $%
\epsilon $ is sufficiently small, we get inequality (6.1).
\end{proof}

\begin{proof}[Proof of Lemma 2]
This lemma is proved in [\textbf{\ }Masuda.K and\ Takahashi. K\textbf{\ }%
\cite{Masuda}, Lemma 2.2].
\end{proof}

\begin{proof}[Proof of Lemma 3]
Immediate from the maximum principle.
\end{proof}


\begin{thebibliography}{99}
\bibitem{abdelmalek1} \textbf{\ Abdelmalek. S }and\textbf{\ Kouachi. S}, 
\textit{A Simple Proof of Sylvester's (Determinants) Identity}, App.Math.
scie. Vol. 2.2008. no 32. p 1571-1580.

\bibitem{Laurent} \textbf{\ Desvillettes. L and Fellner. K}, \textit{Entropy
Methods for Reaction-Diffusion Systems: Degenrate Diffution}, Discrete and
Continuous Dynamical Systems, Supplement Volume 2007.

\bibitem{Friedman} \textbf{\ Friedman. A}, \textit{Partial Differential
Equations of Parabolic Type}. Prentice Hall Englewood Chiffs. N. J. 1964.

\bibitem{Gierer} \textbf{\ Gierer. A }and\textbf{\ Meinhardt. H}, A Theory
of Biological Pattern Formation. Kybernetik, 1972,12:30-39.

\bibitem{Haraux} \textbf{\ Haraux. A }and\textbf{\ Kirane. M}, \textit{%
Estimations C}$^{1}$\textit{pour des probl\`{e}mes paraboliques semi-lin\'{e}%
aires}, Ann. Fac. Sci. Toulouse 5 (1983), 265-280.

\bibitem{Henry} \textbf{\ Henry. D}, Geometric \textit{Theory of Semi-linear
Parabolic Equations}. Lecture Notes in Mathematics 840, Springer-Verlag,
New-York, 1984.

\bibitem{jiang} \textbf{\ jiang. H}, \textit{Global existence of \ }Solution 
\textit{of \ an }Activator-Inhibitor System, Discrete and continuous
Dynamical \ Systems. V14,N4 April 2006.p 737-751.

\bibitem{Wu} \textbf{\ Jianhua. W }and\textbf{\ Yanling. L}, \textit{Global
Classical Solution for the Activator-Inhibitor Model}. Acta Mathematicae
Applicatae Sinica (in Chinese), 1990, 13: 501-505.

\bibitem{Masuda} \textbf{\ Masuda.K }and\textbf{\ Takahashi. K}, \textit{%
Reaction-diffusion systems in the Gierer-Meinhardt theory of biological
pattern formation}. Japan J. Appl. Math., 4(1): 47-58, 1987.

\bibitem{Meinhardt3} \textbf{Meinhardt. H, Koch. A }and\textbf{\ Bernasconi.
G}, \textit{Models of pattern formation applied to plant development},
Reprint of a chapter that appeared in: Symmetry in Plants (D. Barabe and R.
V. Jean, Eds), World Scientific Publishing, Singapore;pp. 723-75.

\bibitem{Mingde} \textbf{\ Mingde. L, Shaohua. C }and\textbf{\ Yuchun. Q}, 
\textit{Boundedness and Blow Up for the general Activator-Inhibitor Model,}
Acta Mathematicae Applicatae Sinica, vol.11 No.1. Jan, 1995.

\bibitem{W. Ni} \textbf{Ni. W, Suzuki. K }and\textbf{\ Takagi. I}. \textit{%
The dynamics of a kinetic activator--inhibitor system}. J. Differential
Equations 229 (2006) 426--465.

\bibitem{Pazy} \textbf{\ Pazy. A}, \textit{Semigroups of Linear Operators
and Applications to Partial Differential Equations}, Applied Math. Sciences
44, Springer-Verlag, New York (1983).

\bibitem{Rothe} \textbf{\ Rothe}. \textbf{F.}\textit{Global Solutions of
Reaction-Diffusion Equations}. Lectur Notes in Mathematics ,1072,
Springer-Verlag, Berlin, 1984.

\bibitem{Smoller} \textbf{\ Smoller. J}, \textit{Shock Waves and
Reaction-Diffusion Equations}, Springer-Verlag, New York (1983).

\bibitem{Abraham} \textbf{Trembley. A}, \textit{Memoires pour servir a
l'histoire d'un genre de polypes d'eau douce}, abras en forme de cornes.
1744.

\bibitem{Turing} \textbf{\ Turing. A. M}, \textit{The chemical basis of
morphogenesis}. Philosophical Transactions of the Royal Society (B), 237:
37-72, 1952.
\end{thebibliography}
\end{document}